\documentclass[12pt]{amsart} 

\newif\ifpdf
\pdffalse
\ifx\pdfoutput\undefined
\else
  \ifx\pdfoutput\relax
    \let\pdfoutput\undefined
  \else
    \ifcase\pdfoutput
    \else
      \pdftrue
    \fi
  \fi
\fi

\ifpdf
  \usepackage[pdftex]{graphicx}
\else
  \usepackage[dvips]{graphicx}
\fi

\ifpdf
 \DeclareGraphicsExtensions{.pdf, .png, .jpg, .jpeg, .tif, .tiff}
\else
 \DeclareGraphicsExtensions{.eps,.ps}
\fi

\usepackage[T1]{fontenc} 
\usepackage[latin1]{inputenc} 
\usepackage{amsmath} 
\usepackage{amsthm} 
\usepackage{amsfonts} 
\usepackage{bbm} 
\usepackage[english]{babel} 
\usepackage{varioref} 
\usepackage{enumerate} 
\usepackage{a4wide}

\usepackage{ae}

\theoremstyle{definition}

\theoremstyle{remark}

\theoremstyle{plain}
\newtheorem{thm}{Theorem}

\newcommand{\norm}[1]{\ensuremath{\left\Vert #1 \right\Vert}}
\newcommand{\abs}[1]{\ensuremath{\left\vert #1 \right\vert}}

\newcommand{\natnum}{\ensuremath{\mathbb{N}}}
\newcommand{\integer}{\ensuremath{\mathbb{Z}}}

\newcommand{\reals}{\ensuremath{\mathbb{R}}}

\newcommand{\x}{\mathbf{x}}
\newcommand{\p}{\mathbf{p}}
\DeclareMathOperator{\dist}{dist}

\title[Diophantine approximation and planar distance functions]{Metric
  Diophantine approximation with respect to planar distance functions}
\author[Simon Kristensen]{Simon Kristensen (Edinburgh)}
\address{School of Mathematics, University of Edinburgh, King's
  Buildings, JCMB, Mayfield Road, Edinburgh, EH9 3JZ, Scotland}
\email{Simon.Kristensen@ed.ac.uk}
\thanks{The author is a William Gordon Seggie Brown Fellow}
\date{\today}


\begin{document}

\begin{abstract}
  We outline a proof of an analogue of Khintchine's Theorem in
  $\mathbb{R}^2$, where the ordinary height is replaced by a distance
  function satisfying an irrationality condition as well as certain
  decay and symmetry conditions.
\end{abstract}

\maketitle

\section{Introduction}
\label{sec:introduction}

A star body in the plane $S \subseteq \reals^2$ is defined as an open
set which has the origin as an interior point, such that for any $\x
\in \reals^2$, there exists a $t_0 \in (0, \infty]$ such that for $t <
t_0$, $t \x \in S$ and for $t > t_0$, $t \x \notin S$ (see
\cite{MR28:1175}). To such sets, one can associate a distance function
$F: \reals^2 \rightarrow \reals$ such that $F(\x) \geq 0$ for all $\x
\in \reals^2$, and such that for any $\x \in \reals^2$ and any $t \geq
0$, $F(t \x) = t F(\x)$. The open star body $S$ may be expressed as
the set of points satisfying $F(\x) < 1$. Conversely, to each
continuous distance function, we may associate an open star body
through this relation. Throughout this paper, vectors will be denoted
in bold typeface and their coordinates in normal typeface with a
subscript, \emph{e.g.}, $\x = (x_1, x_2) \in \reals^2$.

Let $F: \mathbb{R}^2 \rightarrow \mathbb{R}$ be a distance function
and let $\psi: \mathbb{N} \rightarrow \mathbb{R}_+$ be a function such 
that $q \psi(q)$ is non-increasing. The measure of the set
\begin{multline}
  \label{eq:1}
  W(F;\psi) = \big\{\x \in [0,1)^2 :  \text{ for
    infinitely many } q \in \natnum, \\
  \text{ there is a } \p \in \integer^2, \text { such that }
    F(\x-\p/q) < \psi(q) \big\}, 
\end{multline}
has been studied for a number of distance functions $F$, the most
important of which are the height $\max\{\abs{x_1}, \abs{x_2}\}$ and
the multiplicative distance function $\sqrt{\abs{x_1}\abs{x_2}}$. For
both of these, the measure is null or full accordingly as a certain
series depending on $F$ and $\psi$ converges or diverges. For the
height, the series is $\sum_{q=1}^\infty q^2 \psi(q)^2$. This was
shown by Khintchine
\cite{khintchine24:_einig_saetz_ketten_anven_theor_dioph_approx}.  For
the multiplicative distance function, the series is $\sum_{q=1}^\infty
q^2 \psi(q)^2 \log(q^{-1} \psi(q)^{-1})$. This result is a corollary
of a result by Gallagher \cite{MR28:1167}.

Note that in the original statements of theorems of Khintchine and
Gallagher, the sets are defined slightly differently. Originally, $\p$
is the integral vector closest to $q \x$. However, using the form of
the distance functions, it may easily be shown that the sets
originally studied by Khintchine and Gallagher are equal to
$W(F;\psi)$ for the appropriate distance function.

The Hausdorff dimension of the set $W(F;\psi)$ was studied by Dodson
and Hasan \cite{MR92h:11063, MR93g:11078} under a technical covering
condition. When the dimension is less than $2$, the set is obviously
null. Unfortunately, this covering condition is very difficult to
check for most non-trivial star bodies. Also, their results do not
provide conditions under which the measure of the set is positive.

We will outline a proof of the fact that for a class of distance
functions, the measure of $W(F;\psi)$ is full regardless of the error
function $\psi$. We will make a number of technical assumptions on the
distance function, some of which are unnecessary in general, but which
will make the present presentation simpler. The full result will be
published elsewhere in a joint paper with M.~M.~Dodson \cite{dodson}
where additional distance functions will also be treated.

Let $F: \mathbb{R}^2 \rightarrow \mathbb{R}_+$ be a distance function
satisfying the following conditions:
\begin{enumerate}
\item \label{item:1} The set $F^{-1}(0)$ consists of a family of
  half-lines such that at least one has an irrational slope.
\item \label{item:2} One irrationally sloped half-line $L$ from
  $F^{-1}(0)$ satisfies that for any $M > 0$,
  \begin{displaymath}
    \mu\left(\left\{\x \in \reals^2 : F(\x) < 1\right\} \cap \left\{\x 
        \in \reals^2 : \dist(\x,L) < M\right\} \right) = \infty.
  \end{displaymath}
\item \label{item:3} For any point $\x$ on the half-line $L$ with
  $\abs{\x}$ sufficiently large, the line orthogonal to $L$ through
  $\x$ intersects the level curve $F(\x) = 1$ in such a way that the
  two intersection points closest to $\x$ are at equal distance from
  $\x$.
\item \label{item:4} The distance from $L$ to the closest intersection
  point of the line orthogonal to $L$ through the point $\x$ is
  decreasing as a function of $\abs{\x}$ for $\abs{\x}$ large enough.
\end{enumerate}
Broadly speaking, the conditions mean that the star body has
\eqref{item:1} an unbounded component about an irrationally sloped
half line, \eqref{item:2} such that this component has infinite
measure, \eqref{item:3} such that the unbounded component is symmetric
about the irrationally sloped line and \eqref{item:4} such that the
width of the star body decreases as we move along the unbounded
component. Condition \eqref{item:3} is strictly speaking not needed,
but it will make the present presentation easier to follow. The full
result without the additional requirements will be proved in
\cite{dodson}. 

We now state the main result.
\begin{thm}
  \label{thm:main-theorem}
  Let $F: \reals^2 \rightarrow \reals$ be a distance function
  satisfying conditions \eqref{item:1}--\eqref{item:4}. Then, for any
  $\psi:\mathbb{N} \rightarrow \mathbb{R}_+$, for any $q \in \natnum$
  and for almost all $\x \in [0,1)^2$ there are infinitely many $\p
  \in \integer^2$ such that
  \begin{equation}
    \label{eq:3}
    F(\x - \p/q) < \psi(q).
  \end{equation}
  In particular, 
  \begin{displaymath}
    \abs{W(F;\psi)} = 1.
  \end{displaymath}
\end{thm}
This result extends Khintchine's Theorem in a rather surprising way to
a larger class of distance functions not covered by any of the
previous results. Gallagher's result \cite{MR28:1167} allows us to
deal with distance functions where $F^{-1}(0)$ may contain any of the
coordinate half-axes but no other half-lines. Furthermore, as the
measure is always full, the results of Dodson and Hasan are not
applicable to these situations.

\section{Outline of proof}
\label{sec:outline-proof}

We now outline the proof of Theorem \ref{thm:main-theorem}. Let $q \in
\mathbb{N}$ be fixed and let $\alpha$ denote the slope of the line $L$
satisfying assumptions \eqref{item:1}--\eqref{item:4}. Suppose that
$\alpha \geq 1/2$ and that $L$ is in the first quadrant. This causes
no loss of generality, as we may change the role of the axes if this
is not the case to obtain a similar situation. The family of lines
$\mathcal{F} = \{L+\p/q : \p \in \mathbb{Z}^2\}$ induces a family of
dense geodesics on the standard torus $\mathbb{T}^2$ having $[0,1)^2$
as its fundamental domain. This follows from Kronecker's Theorem (see
\emph{e.g.}  \cite[Proposition 1.5.1]{MR96c:58055}).

We fix one such geodesic and consider a fixed horizontal line across
$[0,1)^2$. With the endpoints identified, this may be identified with
a unit circle $\mathbb{S}$. Furthermore, the geodesic on $\mathbb{T}^2$
intersects this circle at multiples of $\alpha^{-1}$ plus some
constant corresponding to the first intersection $y_0$. That is, we
have points $\{y_0 + n \alpha^{-1}\}$, where $\{x \}$ denotes the
fractional part of $x$. These points are uniformly distributed, but
this is unfortunately not sufficient to ensure that the measure of the
intersection between $W(F;\psi)$ and the circle is full. 

We now study the intersection of the unbounded branch of the star
bodies centred at the points $\p/q$ and the circle. About the geodesic,
there is a strip which becomes increasingly narrow as we move along.
The width of this strip is decreasing by assumption \eqref{item:4}.
Considering the intersection with $\mathbb{S}$, we find that to each
point $\{y_0 + n \alpha^{-1}\}$ corresponds an interval of some
length. We consider subintervals of these centred at $\{y_0 + n
\alpha^{-1}\}$ of width $\rho_n$ comparable to the width of the strip
at a point along $L$ at some multiple of $n$. That $\rho_n$ is
comparable to the width follows from assumption \eqref{item:3}.
Finally, assumption \eqref{item:2} along with the homogeneity of the
distance function implies that $\sum_{n=1}^{\infty} \rho_n = \infty$
regardless of the error function $\psi$.

We have now reduced the problem of finding the measure of the set of
points $\x$ on the horizontal line chosen for which $F(\x-\p/q) <
\psi(q)$ for infinitely many $\p \in \mathbb{Z}^2$ to the problem of
finding the measure of the set of points $x \in \mathbb{S}$ satisfying
the inequalities
\begin{equation}
  \label{eq:2}
  \norm{n \alpha^{-1} + y_0 - x} < \rho_n
\end{equation}
for infinitely many $n \in \natnum$. Indeed, such points will satisfy
the original inequalities by the above arguments. 

When the series $\sum \rho_n$ converges, the Hausdorff dimension of
the set of points $y_0 - x$ satisfying \eqref{eq:2} infinitely often
was found by Bugeaud \cite{MR1972699} using regular systems. His
arguments may be adapted to deal with the situation when the series is
divergent. For this purpose, the notion of a locally ubiquitous system
is required.

Let $\lambda:\reals_+ \rightarrow \reals_+$ be a function decreasing
to zero and let $(N_r)$ be a strictly increasing sequence of natural
numbers. The set of points $\mathcal{R} = \{z_n: n \in \natnum\}$ in
the metric space $\mathcal{S}$ is said to be \emph{locally ubiquitous
  relative to $\lambda$ and $(N_r)$} if there is a constant $\kappa >
0$ such that for any interval $I=(c-\rho, c+\rho)$ with $\rho$
sufficiently small and any $r \geq r_0 (I)$,
\begin{displaymath}
  \mu\left(I \cap \bigcup_{N_r \leq n < N_{r+1}} \left(z_n -
      \lambda(N_r), z_n + \lambda(N_r)\right) \right) \geq  2 \kappa
  \rho. 
\end{displaymath}
This is a quantitative version of the local density of the points in
$\mathcal{R}$. 

Using the so-called Three Distances Theorem due to S\'os
\cite{MR20:34}, we may easily prove that the existence of a strictly
increasing sequence $(N_r)$ such that the set of points $\{n
\alpha^{-1} + y_0\}$ is locally ubiquitous with respect to $\lambda$
and $(N_r)$, where $\lambda(N) = 3/(1+N_r)$ for $N \in [N_r,
N_{r+1})$. Using a recent result due to Beresnevich, Dickinson and
Velani \cite{beresnevich03:_measur}, this along with the divergence of
$\sum \rho_n$ is sufficient to ensure that the set of $x \in
\mathbb{S}$ for which \eqref{eq:2} is satisfied for infinitely many $n
\in \mathbb{N}$ has full measure. 

As we may prove full measure for every horizontal line across $[0,1)$,
the theorem follows on integrating over the relevant characteristic
functions in the vertical direction.

\section{Concluding remarks}
\label{sec:concluding-remarks}

We have only given an outline of the proof of the main theorem. In
fact, a more general theorem may be proven, where assumption
\eqref{item:3} is no longer needed. In \cite{dodson} we will prove the
more general result in full detail. The ubiquity argument becomes more
involved when \eqref{item:3} is not assumed.

In this paper, we will also prove a Khintchine type theorem for
distance functions where $F^{-1}(0)$ consists of finitely many
half-lines with rational slopes. In this case, the result depends on
the error function and takes a form resembling the classical theorem.
This shows that the phenomenon studied by the various Khintchine type
theorems is very sensitive to minor changes in the distance function.
For example, a small rotation of the star body corresponding to the
multiplicative distance function may cause a slope to go from rational
to irrational, and thus completely changing the nature of the $0$--$1$
law of the associated Khintchine type theorem. This is another example
of the wide contrast between the rational and the irrational.

\section{Acknowledgements}
\label{sec:acknowledgements}

I thank Victor Beresnevich, Maurice Dodson and Sanju Velani for their
helpful comments and suggestions. 

\providecommand{\bysame}{\leavevmode\hbox to3em{\hrulefill}\thinspace}
\providecommand{\href}[2]{#2}

\end{document}